\documentclass[12pt]{article}

\usepackage{amsmath,amssymb,amsthm}

\newtheorem{theorem}{Theorem}[section]
\newtheorem{lemma}[theorem]{Lemma}

\newtheorem{corollary}[theorem]{Corollary}

\theoremstyle{definition}

\newtheorem{definition}{Definition}[section]

\newcommand{\ba}{\mathbf{a}}
\newcommand{\bb}{\mathbf{b}}
\newcommand{\bz}{\mathbf{z}}
\newcommand{\bbf}{\mathbf{f}}

\newcommand{\bA}{\mathbf{A}}
\newcommand{\bB}{\mathbf{B}}
\newcommand{\bC}{\mathbf{C}}
\newcommand{\bE}{\mathbf{E}}

\newcommand{\bI}{\mathbf{I}}

\newcommand{\bF}{\mathbf{F}}
\newcommand{\fF}{\mathfrak{F}}
\newcommand{\fP}{\mathfrak{P}}

\newcommand{\bR}{\mathbb{R}}
\newcommand{\bS}{\mathbb{S}}
\newcommand{\bT}{\mathbb{T}}
\newcommand{\bZ}{\mathbb{Z}}

\newcommand{\cF}{\mathcal{F}}
\newcommand{\cH}{\mathcal{H}}

\newcommand{\cL}{\mathcal{L}}

\newcommand{\cP}{\mathcal{P}}

\newcommand{\cS}{\mathcal{S}}

\newcommand{\nn}{\nonumber}

\newcommand{\hb}{\hbar}
\newcommand{\pa}{\partial}

\newcommand{\ve}{\varepsilon}
\newcommand{\ol}{\overline}

\newcommand{\opa}{\overline{\partial}}
\newcommand{\od}{\overset{\rm def}{=}}

\renewcommand{\Re}{\operatorname{Re}}

\renewcommand{\div}{\operatorname{div}}

\newcommand{\ad}{\operatorname{ad}}
\newcommand{\su}{\operatorname{su}}
\newcommand{\so}{\operatorname{so}}
\newcommand{\tr}{\operatorname{tr}}
\newcommand{\const}{\operatorname{const}}

\def\zs#1{_{\lower2pt\hbox{$\scriptstyle#1$}}}

\begin{document}

\title{Resonance Gyrons and Quantum Geometry}

\author{Mikhail Karasev\thanks{This work was partially supported
by RFBR (grant 05-01-00918-a).}\\ \\
\it Moscow Institute of Electronics and Mathematics\\
karasev@miem.edu.ru \\ \\ \\
\hfill \it Dedicated to Hideki Omori}

\date{}
                                                
\maketitle

\begin{abstract}
We describe irreducible representations, coherent states and
star-products for algebras of integrals of motions (symmetries)
of two-dimensional resonance oscillators. We demonstrate how
the quantum geometry (quantum K\"ahler form, metric, quantum
Ricci form, quantum reproducing measure) arises in this problem.
We specifically study the distinction between the isotropic
resonance $1:1$ and the general $l:m$ resonance 
for arbitrary coprime $l,m$.
Quantum gyron is a dynamical system in the resonance algebra. 
We derive its Hamiltonian in irreducible representations and
calculate the semiclassical asymptotics of the gyron spectrum
via the quantum geometrical objects.
\end{abstract}


\bigskip

\section{Introduction}

For complicated dynamical systems, it is important to be able 
to abstract from studying concrete motions or states and 
to observe surrounding structures, like spaces, algebras, etc., 
which carry essential properties of the variety of motions 
in the whole.

For quantum (wave) systems, the standard accompanying
mathematical structures are algebras of ``observables,''
i.e., functions on phase spaces, 
and representations of these algebras in Hilbert vector spaces
of ``states.'' 
This is the starting viewpoint for the mathematical quantization
theory \cite{j1}--\cite{j10}.
The more complicated systems are studied the more complicated
algebras and phase spaces (symplectic manifolds) have to be
used. Note that for general symplectic and even K\"ahlerian
manifolds the quantization problem is still unsolved.

It was demonstrated in \cite{j11,j12} that for general
symplectic manifolds it is possible to approximate 
the symplectic potential by its quadratic part 
(the oscillator!), then to use this quadratic part 
in order to define the standard Groenewold--Moyal \cite{j13,j14}
product on the tangent spaces, and to construct a formal
$*$-product on the original manifold by a perturbation theory. 
Such oscillator-generated quantum manifolds were called 
the ``Weyl manifolds'' in~\cite{j11}.

In quantum and wave mechanics, one often meets 
a situation similar in certain sense: 
the dynamics of a system is, in general, chaotic, 
but there are some exclusive invariant submanifolds  
(for instance, equilibrium points) in the phase space around
which the dynamics is regular and can be approximated by the
oscillator motion in directions transversal to the
submanifold. 
Thus the given system contains inside a built-in harmonic
oscillator plus certain anharmonic part near the equilibrium:
\begin{equation} 
\frac12\sum(p^2_j+\omega^2_j q^2_j)+\text{cubic}+\text{quartic}+\dots
\tag{1.1}
\end{equation}

If the frequencies $\omega_j$ of the harmonic part are
incommensurable (not in a resonance), then in a small
neighborhood of the submanifold the anharmonic part just
slightly perturbs these frequencies, and the whole motion is
performed along the perturbed Liouville tori. This is the
well-investigated situation both on the classical and quantum
levels \cite{j15}--\cite{j21}.

If the frequencies $\omega_j$ are in a {\it resonance} then all
standard approaches do not work and the picture occurs to be
much more interesting from the viewpoint of quantum geometry.
Here we will follow the works \cite{j22}--\cite{j25}.

First of all, in the resonance case the Liouville tori are
collapsed (to a smaller dimension), and the anharmonic part
generates a nontrivial ``averaged'' motion in the new phase
spaces: in the symplectic leave $\Omega$ of the commutant 
$\cF_\omega$ of the harmonic part.
The new phase spaces represent certain hidden dynamics 
committed to the resonance. This dynamics describes a {\it
precession} 
of the parameters of the resonance harmonic motion under the
action of the anharmonic part. We call this dynamical system a
{\it gyron} (from the Greek word ``gyro,'' i.e., ``rotating'').

In the simplest case of the isotropic 1:1 resonance for two
degrees of freedom the gyron system is just the Euler top system
from the theory of rigid body rotations, which is related to the
linear Poisson brackets. For the general $l:m$ resonance, 
{\it the gyron is described by a nonlinear Poisson brackets
polynomial of degree} $l+m-1$, see in \cite{j24,j25}.

Of course, in the quantum case the resonance function algebra 
$\cF_\omega$ has to be replaced by a resonance operator algebra
$\fF_\omega$ which consists of operators commuting with the
quantum oscillator $\frac12\sum_j(\hat p^2_j+\omega^2_jq^2_j)$,
where $\hat p_j=-i\hb\pa/\pa q_j$.  
This algebra is described by nonlinear commutation relations
of polynomial type, see in \cite{j24,j25}. 
It is the dynamic algebra for {\it quantum gyrons}.

Note that there is a variety of important physical models
containing inside the resonance Hamiltonians like (1.1). The
quantum gyrons in these models can be considered as an analog of
known quasiparticles similar to polarons, rotons, excitons, 
etc.\footnote{Attention to this was paid by V.~Maslov.}
As the simplest example, we mention the models of nano-physics
(quantum dots, artificial atoms, quantum wires, 
see examples in \cite{j24}). 
Another example is the fiber waveguides in optics; they are
described by the Hamiltonian
\begin{equation}
p^2-n^2(q),\qquad q,p\in\bR^3,
\tag{1.2}
\end{equation}
where $n(q)$ is the refraction index having the maximum value 
along the waveguide axis, that is, along an arbitrary smooth
curve in $\bR^3$. 
The quadratic part of $n^2(q)$ in directions transversal to
this curve is assumed to have commensurable frequencies in a
certain resonance proportion $\omega_1:\omega_2=l:m$, where
$l,m$ are coprime integers. The quantum gyron in this model
describes certain hidden ``polarization'' of the light beam
along the given curve in the optical medium, see in \cite{j24}.
The propagation of such {\it optical gyrons} and their spectrum 
depend on the anharmonic part of the refraction index, and so
one can control the properties of the gyron waves by changing
the geometry of the curve just by bending the optical fiber. 

The aim of the given paper is to describe the quantum geometry
of the gyron phase spaces in the case of the $l:m$ resonance. 

If $l=m=1$, then these phase spaces $\Omega$ are just
homogeneous spheres $\bS^2$, that is, the coadjoint $su(2)$
orbits.  
The quantum geometry in this case coincides with the classical
symplectic (K\"ahlerian) geometry generated by linear
Lie--Poisson brackets. 

If at least one of the integers $l$ or $m$ exceeds~$1$, then,
as we will see below, the quantum geometry occurs to be unusual.
The quantum phase spaces are still diffeomorphic to $\bS^2$,
but the classical symplectic form is singular on them. 
The correct symplectic (K\"ahlerian) form and the reproducing
measure of the quantum phase space are chosen from the
nontrivial condition that the operators of irreducible
representations of the quantum resonance algebras
$\fF_\omega=\fF_{l,m}$ have to be differential operators, 
not pseudodifferential (the maximal order of these operators is
$\max(l,m)$). 

Thus the geometry \cite{j26,j27} determining the
Wick--Klauder--Berezin 
$*$-product on the gyron phase space has a purely quantum
behavior and the $*$-product itself cannot be obtained by a
formal deformation technique from the classical data.

Note that here we mean the phase spaces corresponding 
either to low energy levels of the oscillator 
(i.e., to the nano-zone near its equilibrium point, 
in the terminology of \cite{j24}) 
or to excited levels 
(i.e., to the micro-zone).
Thus one can talk about {\it quantum nano- or micro-geometry}
generated by the $l:m$ frequency resonance.

The distinction between the specific case $l=m=1$ 
and the generic case $\max(l,m)>1$ is the distinction between
algebras with linear and nonlinear commutation relations. 
We see that the nonlinearity of relations in the algebra
$\fF_{l,m}$ (the absence of a Lie group of symmetries) 
for the resonance oscillator implies the quantum character 
of the phase spaces in nano- and micro-zones near the ground
state. The motion in these spaces is the gyron dynamics. 
In the nano-zone, this dynamics is purely quantum and does not
have a classical analog at all.
In the micro-zone, the gyron dynamics and the gyron spectrum
can be described by semiclassical methods \cite{j23,j24} if one
at first fixes the quantum geometry of the gyron phase space.

Applying this theory, for instance, to optical gyrons, we come
to the conclusion that the {\it light beam propagating near the
axis of a resonance fiber waveguide cannot be described by
purely geometric optics and carry essentially quantum properties}.
This opens an opportunity to apply such simple optical devices,
for example, in constructing elements of quantum computers.

Also note that the $l:m$ resonance oscillators, which we discuss
here, can be presented in the form 
\begin{equation}
\hat{l}+\hat{m},
\tag{1.3}
\end{equation}
where $\hat{l}$ and $\hat{m}$ are mutually commutating action
operators with spectra $l\cdot\bZ_+$ and $m\cdot\bZ_+$ 
in the Hilbert space $\cL=L^2(\bR\times\bR)$.
The operators $\hat{l}$ and $\hat{m}$ can be considered as
``quantum integer numbers'' and their sum as a quantum sum of
integers. Then the representation theory of the algebra
$\fF_{l,m}$ and the corresponding quantum geometry could be
considered as a brick to construction of something like 
``quantum arithmetics.'' 

\section{Commutation relations and Poisson brackets 
for $l:m$ resonance}

The Hamiltonian of the resonance oscillator (1.3) 
can be written as 
\begin{equation}
\bE=l \bb^*_1 \bb_1 +m \bb^*_2 \bb_2.
\tag{2.1}
\end{equation}
Here $l,m$ are coprime integers, 
$\bb_1,\bb_2$ are annihilation operators
in the Hilbert space $\cL$, 
and $\bb^*_1,\bb^*_2$ are the conjugate
creation operators. The commutation relations are 
$$
[\bb_1,\bb^*_1]=[\bb_2,\bb^*_2]=\hb \bI,
$$
all other commutators are zero. 

In the algebra generated by 
$\bb_1$, $\bb_2$, $\bb^*_1$, $\bb^*_2$, 
let us consider the commutant of the element (2.1).
This commutant is a nontrivial, noncommutative subalgebra. 
We call it a {\it resonance algebra}. 
It is related to quantum gyrons.

Note that the resonance algebra is generated by the following
four elements:
\begin{equation}
\bA_1=\bb^*_1\bb_1,\qquad
\bA_2=\bb^*_2\bb_2,\qquad
\bA_+=(\bb^*_2)^l\bb^m_1,\qquad
\bA_-=\bA^*_+.
\tag{2.2}
\end{equation}

Let us define the following polynomials 
\begin{align}
\rho(A_1,A_2)&\od \prod^{m}_{j=1}(A_1+j\hb)\cdot 
\prod^{l}_{s=1}(A_2-s\hb+\hb),
\tag{2.3}\\
\varkappa(A_1,A_2)&\od lA_1+mA_2.
\nn
\end{align}

\begin{lemma}
Elements {\rm(2.2)} obey the commutation relations
\begin{gather}
[\bA_1,\bA_2]=0
\nn\\
[\bA_1,\bA_{\pm}]=\mp\hb m \bA_{\pm},\qquad
[\bA_2,\bA_{\pm}]=\pm\hb l \bA_{\pm},
\tag{2.4}\\
[\bA_-,\bA_+]=\rho(\bA_1-\hb m,\bA_2+\hb l)-\rho(\bA_1,\bA_2).
\nn
\end{gather}
\end{lemma}

\begin{lemma}
In the abstract algebra $\fF_{l,m}$ with relations {\rm(2.4)}
there are two Casimir elements 
$$
\boldsymbol\kappa=\varkappa(\bA_1,\bA_2),\qquad 
\bC=\bA_+\bA_- -\rho(\bA_1,\bA_2).
$$
In realization {\rm(2.2)} the Casimir element $\bC$ is
identically zero, and the Casimir element $\boldsymbol\kappa$
coincides with the oscillator Hamiltonian $\bE$ {\rm(2.1)}. 
\end{lemma}

Note that the operators $\bA_1,\bA_2$ (2.2) are self-adjoint,
but $\bA_+$ is not. Let us introduce the self-adjoint operators 
$\bA_3,\bA_4$ by means of the equalities
$$
\bA_{\pm}=\bA_3\mp i\bA_4.
$$
Then commutation relations (2.4) read
\begin{align}
[\bA_1,\bA_2]&=0,\qquad [\bA_1,\bA_3]=i\hb m\bA_4,\qquad
[\bA_1,\bA_4]=-i\hb m\bA_3,
\nn\\
&\qquad\qquad\,\,  [\bA_2,\bA_3]=-i\hb l\bA_4,\quad\,\,\,\,
[\bA_2,\bA_4]=i\hb l\bA_3,
\tag{2.4a}\\
[\bA_3,\bA_4]&
=\frac{i}2\big(\rho(\bA_1-\hb m,\bA_2+\hb l)-\rho(\bA_1,\bA_2)\big).
\nn
\end{align}

Let us denote by $A_j$ the classical variable 
(a coordinate on $\bR^4$) corresponding to
the quantum operator $\bA_j$. Then the relations (2.4a) are
reduced to the following Poisson brackets on $\bR^4$:
\begin{align}
\{A_1,A_2\}&=0, 
\nn\\
\{A_1,A_3\}&=-m A_4,\qquad \{A_1,A_4\}=m A_3,
\tag{2.5}\\
\{A_2,A_3\}&=l A_4,\qquad\quad\,\, \{A_2,A_4\}=-l A_3,
\nn\\
\{A_4,A_3\}&=\frac12(l^2 A_1-m^2A_2) A^{m-1}_1 A^{l-1}_2.
\nn
\end{align}

\begin{lemma}
Relations {\rm(2.5)} determine the Poisson brackets on $\bR^4$
with the Casimir functions
$$
\varkappa=lA_1+mA_2,\qquad 
C=A^2_3 + A^2_4 - A^m_1A^l_2.
$$
\end{lemma}

\begin{lemma}
In the subset in $\bR^4$ determined by the inequalities 
$A_1\geq0$ and $A_2\geq0$, there is a family of surfaces
\begin{equation}
{\Omega}=\{\varkappa=E,C=0\},\qquad E>0,
\tag{2.6}
\end{equation}
which coincide with the closure of symplectic leaves $\Omega_0$ 
of the Poisson structure {\rm(2.5)}.
These surfaces are diffeomorphic to the sphere{\rm:} 
${\Omega}\approx\bS^2$.

The topology of the symplectic leaves $\Omega_0$ is the
following{\rm:} 
\begin{itemize}
\item[--] 
if $l=m=1$, then $\Omega_0=\Omega${\rm;}
\item[--] 
if $l=1$, $m>1$ or $l>1$, $m=1$, then $\Omega_0$ is obtained
from $\Omega$ by deleting the point 
$(0,\frac{E}{m},0,0)$ or the point $(\frac{E}{l},0,0,0)${\rm;} 
\item[--] 
if $l>1$, $m>1$, then $\Omega_0$ is obtained from $\Omega$ 
by deleting both the points $(0,\frac{E}{m},0,0)$ and 
$(\frac{E}{l},0,0,0)$.
\end{itemize}
\end{lemma}

\begin{lemma}
If $l>1$ or $m>1$, then the Kirillov symplectic form
$\omega_0$ on the leaf $\Omega_0\subset\Omega$ has a weak
{\rm(}integrable{\rm)} singularity at the point $A_2=0$ 
or $A_1=0$.
The symplectic volume of $\Omega_0$ is finite
\begin{equation}
\frac{1}{2\pi}\int_{\Omega_0}\omega_0=\frac{E}{lm}.
\tag{2.7}
\end{equation}
\end{lemma}

\begin{lemma}
On the subset $A_1>0$ the complex coordinate 
\begin{equation}
z_0=\frac{A_3+i A_4}{A^m_1}
\tag{2.8}
\end{equation}
determines a partial complex structure consistent with the
brackets {\rm{(2.5)}} in the sense of {\rm\cite{j32}}.
On each symplectic leave $\Omega_0$, this partial complex
structure generates the K\"ahlerian structure with the potential 
\begin{equation}
\Phi_0=\int^{|z_0|^2}_{0}\bigg(\frac{E}{2lm}
+\alpha\zs{E}(x)\bigg)\frac{dx}{x},
\qquad \omega_0=i\opa\pa \Phi_0.
\tag{2.9}
\end{equation}
Here $\pa$ is the differential by $z_0$ and  
$\alpha\zs{E}=\alpha\zs{E}(x)$ is the solution of the
equation 
\begin{equation}
x=\bigg(\frac{E}{2m}+l\alpha\zs{E}\bigg)^l
\bigg(\frac{E}{2l}-m\alpha\zs{E}\bigg)^{-m}
\tag{2.10}
\end{equation}
with values on the interval 
$-\frac{E}{2lm}\leq\alpha\zs{E}\leq \frac{E}{2lm}$.

The singular points of $\omega_0$ on $\Omega_0$ correspond to the
poles  
\begin{align}
A_2&=0\quad\Longleftrightarrow\quad z_0=0, \qquad
\omega_0\sim\frac{1}{l^2}\bigg(\frac{E}{l}\bigg)^{m/l}
\frac{dx\wedge d\varphi}{x^{1-1/l}}\quad\text{as}\quad z_0\to0,
\tag{2.11}\\
A_1&=0\quad\Longleftrightarrow\quad z_0=\infty, \qquad
\omega_0\sim\frac{1}{m^2}\bigg(\frac{E}{m}\bigg)^{l/m}
\frac{dx\wedge d\varphi}{x^{1+1/m}}\quad\text{as}\quad z_0\to\infty,
\nn
\end{align}
where $(x,\varphi)$ are polar coordinates, 
$z_0=x^{1/2}\exp\{i\varphi\}$.

The restrictions of coordinate functions to the surface 
{\rm(2.6)} are given by 
\begin{gather}
A_1\bigg|_{\Omega_0}=\frac{E}{2l}-m\alpha\zs{E}(|z_0|^2),
\qquad
A_2\bigg|_{\Omega_0}=\frac{E}{2m}+l\alpha\zs{E}(|z_0|^2),
\tag{2.12}
\\
(A_3+iA_4)\bigg|_{\Omega_0}=z_0\bigg(\frac{E}{2l}
-m\alpha\zs{E}(|z_0|^2)\bigg)^m.
\nn
\end{gather}
\end{lemma}

Note that the properties of classical symplectic
leaves of the $l:m$ resonance algebra, described 
in Lemmas~2.4--2.6, are a particular case of the topology and
geometry of {\it toric varieties}
(in our case the torus $\bT^1=\bS^1$ is the cycle); 
about this see general theorems in \cite{j28}--\cite{j30}. 
The Poisson extension (2.5) by means of polynomial brackets was 
first described in \cite{j22,j23} for the case of $1:2$
resonance and in \cite{j24,j25} for the $l:m$ case, 
as well for the general multidimensional resonances. 
A type of Poisson extension was also considered in \cite{j31}
for some specific class of resonance proportions 
(which does not include, for instance, the $1:2:3$ 
resonance).  

\section{Irreducible representations of $l:m$ resonance algebra}

First of all, let us discuss the basic problems in constructing
irreducible representations of algebras like (2.4), (2.4a). 
Following the standard geometric quantization program \cite{j6} 
one has
to choose a line bundle over symplectic leaves $\Omega_0$ of the
Poisson algebra related to (2.4a), that is,  the Poisson algebra
(2.5). Then this bundle is endowed with the Hermitian connection 
whose curvature is $i\omega_0$, and a Hilbert space $\cH_0$ of
antiholomorphic sections of the bundle is introduced. 
In this Hilbert space,
the operators of irreducible representation of the algebra
(2.4a) are supposed to act and to be self-adjoint.

However, there are two principle difficulties. First, we do not
know which measure on $\Omega_0$ to take in order to determine
the Hilbert norm in the space $\cH_0$. The choice of measure
should imply the {\it reproducing property} \cite{j32,j33}
\begin{equation}
\omega_0=i\opa\pa \ln\sum_k|\varphi^{(k)}_0|^2,
\tag{3.1}
\end{equation}
where $\{\varphi^{(k)}_0\}$ is an
orthonormal basis in $\cH_0$. For the inhomogeneous case, 
where the commutation relations (2.4a) are 
not linear and no Lie group acts on $\Omega_0$, the existence
of such a reproducing measure is, in general, unknown. 
This difficulty was discovered in \cite{j34} (more precisely, it
was observed in \cite{j34} that the Liouville measure generated
by the symplectic form $\omega_0$ does not obey the property
(3.1) in general). 

Secondly, even if one knows the reproducing measure, there is
still a problem: the operators of the irreducible representation
constructed canonically by the geometric quantization scheme
would be pseudodifferential, but not differential operators.
There are additional nontrivial conditions on the complex
structure (polarization) that make the generators of the algebra
be differential operators (of order greater than~$1$, in general).
About such highest analogs of the Blattner--Kostant--Sternberg
conditions for the polarization to be ``invariant'' see in 
\cite{j35,j36}. 

Taking these difficulties into account, we modify the
quantization scheme. From the very beginning, 
we look for an appropriate complex structure and the scalar product
in the space of antiholomorphic functions that guarantee the
existence of an Hermitian representation of the given algebra by
differential operators, and then introduce a ``quantum''
K\"ahlerian form $\omega$ on $\Omega$, a ``quantum'' measure
and the ``quantum'' Hilbert space $\cH$
which automatically obeys the reproducing property
like (3.1) (without ``classical'' label~$0$). 
This approach is explained in \cite{j32,j33,j37}. 

Note that the polynomial structure of the right-hand sides of
relations (2.4), (2.4a) is critically important in this scheme 
to obtain representations by differential operators.

Denote by $\cP_r$ the space of all polynomials 
$\varphi(\lambda)=\sum^{r}_{n=0} \varphi_n\lambda^n$ 
of degree $r\geq0$ with complex coefficients.

\begin{lemma}
Let $f_+$, $f_-$ be two complex functions on $\bZ_+$ 
such that 
\begin{align}
&f_+f_->0\qquad\text{on the subset}\quad\{1,\dots,r\}\subset\bZ_+,
\tag{3.2}\\
&f_-(0)=f_+(r+1)=0.
\nn
\end{align}
Then the differential operators 
\begin{equation}
\ba_+=f_+\bigg(\lambda\frac{d}{d\lambda}\bigg)\cdot \lambda,
\qquad 
\ba_-=\frac1\lambda\cdot f_-\bigg(\lambda\frac{d}{d\lambda}\bigg)
\tag{3.3}
\end{equation}
leave the space $\cP_r$ invariant and they are conjugate 
to each other  with respect to the following scalar product in
$\cP_r${\rm:} 
\begin{equation}
(g,g')\od\sum^{r}_{n=0}\prod^{n}_{s=1}
\frac{\ol{f_-(s)}}{f_+(s)}\varphi_n\ol{\varphi'_n}.
\tag{3.4}
\end{equation}
Any operator $f\big(\lambda\frac{d}{d\lambda}\big)$, 
where $f$ is a real function on $\bZ_+$, is self-adjoint 
in $\cP_r$ with respect to this scalar product. 
\end{lemma}

Now we consider a map
$$
\gamma:\, \bR^k\to \bR^k
$$
and a real function $\rho$ on $\bR^k$. Denote 
by $R_r\subset\bR^k$ the subset of all points
$a_0$ such that 
\begin{align}
\rho(\gamma^{r+1}(a_0))&=\rho(a_0), 
\tag{3.5}\\
\rho(\gamma^{n}(a_0))&>\rho(a_0) \qquad(n=1,\dots,r).
\nn
\end{align}
For any $a_0\in R_r$ we define real functions $f_j$ 
$(j=1,\dots,k)$ on $\bZ_+$ by the formula 
$f_j(n)\od \gamma^n(a_0)_j$, 
and introduce mutually commuting operators in the space $\cP_r$:
\begin{equation}
\ba_j\od f_j\bigg(\lambda\frac{d}{d\lambda}\bigg).
\tag{3.6}
\end{equation}

\begin{lemma}
Let $a_0\in R_r$, and let there be a factorization 
\begin{equation}
\rho(\gamma^n(a_0))-\rho(a_0)=f_+(n)f_-(n),\qquad 
0\leq n\leq r+1,
\tag{3.7}
\end{equation}
where the factors $f_\pm$ obey the property {\rm(3.2)}.
Then the operator $\ba_+$ {\rm(3.3)} and $\ba_j$ {\rm(3.6)} 
in the space $\cP_r$ with the scalar product {\rm(3.4)} satisfy
the relations
$$
\ba^*_+=\ba_-,\qquad \ba^*_j=\ba_j \qquad(j=1,\dots,k),
$$
and 
\begin{align}
&[\ba_j,\ba_s]=0,
\tag{3.8}\\
&\ba_j\ba_+=\ba_+\gamma_j(\ba),\qquad 
\ba_-\ba_j=\gamma_j(\ba)\ba_-\qquad(j=1,\dots,k),
\nn\\
&[\ba_-,\ba_+]=\rho(\gamma(\ba))-\rho(\ba).
\nn
\end{align}
\end{lemma}

\begin{lemma}
Consider the abstract algebra $\fF$ with relations {\rm(3.8)}. 
The element $\bC=\ba_+\ba_--\rho(\ba)$ belongs to the center of
$\fF$. If a function $\varkappa$ on $\bR^k$ is
$\gamma$-invariant, then the element
$\boldsymbol\kappa=\varkappa(\ba)$ belongs to the center of 
$\fF$.

In the representation {\rm(3.3), (3.6)}, these central elements
are scalar: $\bC=\rho(a_0)\cdot\bI$, 
$\boldsymbol\kappa=\varkappa(a_0)\cdot \bI$.
This representation of the algebra $\fF$ is irreducible and 
Hermitian. 

If the map $\gamma$ has no fixed points, then all irreducible
Hermitian representations of the algebra $\fF$ can be obtained
in this way. All such representations of dimension $r+1$ are
parameterized by elements of the set $R_r$ 
{\rm(}$r=0,1,2,\dots${\rm)}.
\end{lemma}

Now let us return to commutation relations (2.4). In this case
$k=2$, the function $\rho$ is given by (2.3), and the mapping 
$\gamma\equiv\Gamma^\hb:\, \bR^2\to\bR^2$ is  
\begin{equation}
\Gamma^\hb\left(\begin{matrix} A_1\\A_2\end{matrix}\right)
\od
\left(\begin{matrix} A_1-\hb m\\ A_2+\hb l \end{matrix}\right).
\tag{3.9}
\end{equation}

It follows from (2.2) that we have to be interested in a subset 
$A_1\geq0$, $A_2\geq0$ in $\bR^2$. 
Also in view of Lemma~2.2, the values of the Casimir element 
$\bC=\rho(a_0)\cdot \bI$ must be zero. From (3.5) we obtain 
\begin{align*}
&\rho(a_0)=\rho(\Gamma^{\hb(r+1)}(a_0))=0,
\\
&\rho(\Gamma^{\hb n}(a_0))>0\qquad (n=1,\dots,r).
\end{align*}
Using (2.3) let us factorize:
\begin{equation}
\rho=\rho_+\rho_-,\qquad
\rho_+(A)\od\prod^{m}_{j=1}(A_1+\hb j),\qquad
\rho_-(A)\od\prod^{l}_{s=1}(A_2-\hb s+\hb).
\tag{3.10}
\end{equation}
It is possible to satisfy (3.7) by choosing 
$$
f_\pm(n)=\rho_\pm(\Gamma^{\hb n}(a_0)).
$$
In this case, the set $R_r\subset \bR^2$ consists of all points 
$a_0=\left(\begin{matrix}\hb(rm+p) \\ \hb q\end{matrix}\right)$ 
for which the pair of integers $p,q$ obeys the inequalities  
\begin{equation}
0\leq q\leq l-1,\qquad 0\leq p\leq m-1.
\tag{3.11}
\end{equation}

The $\gamma$-invariant function $\varkappa$ in our case (3.9) is
just $\varkappa(A)=lA_1+mA_2$. In view of Lemma~3.3, the value of
the second Casimir element $\boldsymbol\kappa=\varkappa(\ba)$ 
in the irreducible representation (3.3), (3.6) is 
$\varkappa(a_0)=E_{r,q,p}$, 
where
\begin{equation}
E_{r,q,p}\od \hb(lmr+lp+mq).
\tag{3.12}
\end{equation}
From Lemma~2.2 we conclude that these numbers coincide 
with eigenvalues of the oscillator $\bE$ (2.1).

Also from (3.4) we see that the scalar product in the space
$\cP_r$ is given by 
\begin{equation}
(\varphi,\varphi')=\sum^{r}_{n=0}\hb^{(l-m)n}
\frac{(q+nl)!(p+(r-n)m)!}{q!(p+rm)!}\varphi_n\ol{\varphi'_n}.
\tag{3.13}
\end{equation}
Thus the vector space of the irreducible representation depends 
on the number $r$ only, but its Hilbert structures are
parameterized by the pairs $q,p$ from (3.11).
That is why below we will use the notation
$\cP_r\equiv\cP_{r,q,p}$. 

Let us summarize the obtained results.

\begin{theorem}
The commutant of the $l:m$ resonance oscillator $\bE$ {\rm(2.1)}
is generated by operators {\rm(2.2)} obeying commutation
relation {\rm(2.4)}. 
The irreducible representation of the algebra {\rm(2.4)},
corresponding to the eigenvalue $E_{r,q,p}$ {\rm(3.12)} of the
operator $\bE$, is given by the following ordinary differential
operators $\ba=(\ba_1,\ba_2)$ and $\ba_\pm${\rm:}
\begin{equation}
\ba=\Gamma^{\hb\lambda \frac{d}{d\lambda}}(a_0),\qquad
\ba_+=\rho_+(\ba)\cdot \lambda,\qquad
\ba_-=\frac1\lambda\cdot\rho_-(\ba).
\tag{3.14}
\end{equation}
Here $a_0=\left(\begin{matrix}\hb(rm+p) \\ \hb q\end{matrix}\right)$, 
the flow $\Gamma$ on $\bR^2$ is defined by {\rm(3.9)} and the
factors $\rho_\pm$ are defined by {\rm(3.10)}. 
The representation {\rm(3.14)} acts in the space $\cP_{r,q,p}$ of 
polynomials in $\lambda$ of degree $r$, 
and it is Hermitian with respect to the scalar product
{\rm(3.13)}. The dimension of this representation is $r+1$.  
\end{theorem}

In fact, formula (3.14) determines just the matrix
representations of the algebra (2.4): 
elements $\ba$ are represented by a diagonal matrix
and $\ba_\pm$ by near-diagonal matrices with respect to the
orthonormal basis of monomials 
\begin{equation}
\varphi^{(k)}(\lambda)=\hb^{(m-l)k/2}
\bigg(\frac{q!(p+rm)!}{(q+kl)!(p+(r-k)m)!}\bigg)^{1/2}\cdot 
\lambda^k\qquad (k=0,\dots,r) 
\tag{3.15}
\end{equation}
in the space $\cP_{r,q,p}$. These matrices are real-valued and
determined by the integer numbers $l,m$ 
(from the resonance proportion) and $r,p,q$
(labeling the representation):
\begin{align}
(\ba_1)_{ns}&=\hb(p+(r-n)m)\delta_{n,s},\qquad 
(\ba_2)_{ns}=\hb(q+nl)\delta_{n,s},
\nn\\
(\ba_+)_{ns}&=\hb^{(l+m)/2}
\bigg(\frac{(q+nl)!(p+(r-s)m)!}{(q+sl)!(p+(r-n)m)!}\bigg)^{1/2}
\delta_{n-1,s},
\tag{3.16}\\
(\ba_-)_{ns}&=(\ba_+)_{sn}.
\nn
\end{align}
Here the matrix indices $n,s$ run over the set $\{0,\dots,r\}$ and 
$\delta_{n,s}$ are the Kronecker symbols.

In the particular case $l=m=1$, from (3.16) one obtains the
well-known Hermitian matrix irreducible representations of the
``spin'' Lie algebra $\su(2)$ with cyclic commutation relation
between generators $\frac12(\bA_1-\bA_2)$, 
$\frac12(\bA_+ + \bA_-)$, $\frac i2(\bA_+ - \bA_-)$.

\section{Quantum geometry of the $l:m$ resonance}

Now we give a geometric interpretation of the obtained
representations of the resonance algebra.

It follows from (3.4) that the element 
$\rho_+(\bA)^{-1}(\bA_3-i\bA_4)$, in the algebra 
generated by relations (2.4),
is represented by the multiplication by $\lambda$ in each
irreducible representation (3.14). If we denote
\begin{equation}
\bz=(\bA_3+i\bA_4)\rho_+(\bA)^{-1},
\tag{4.1}
\end{equation}
then the conjugate operator $\bz^*$ in each irreducible
representation can be taken equal to the multiplication 
by a complex variable $\ol{z}$:
$$
\bz^*=\ol{z}.
$$

Thus, here we change our notation and use $\ol{z}$ instead of
$\lambda$. From now on, $\cP_{r,q,p}$ is the space of anti-holomorphic
functions (polynomials in $\ol{z}$ of degree $\leq r$) on $\bR^2$.

Let us assume that the scalar product (3.13) in the space
$\cP_{r,q,p}$ can be written in the integral form 
\begin{equation}
(\varphi,\varphi')=\frac1{2\pi\hb}\int_{\bR^2}\varphi(\ol{z}(a))
\ol{\varphi'(\ol{z}(a))}L(a)\,da,
\tag{4.2}
\end{equation}
where $da=|d\ol{z}(a)\wedge dz(a)|$ and $a\to z(a)$ 
is the complex coordinate on $\bR^2$. 

\begin{lemma}
The explicit formula for the density $L$ in {\rm(4.2)} 
is 
$$
L(a)=\frac1{4\hb^{rm+p+q+1}(p+rm)!q!x}
\int^\infty_0 A^{rm+p}_1 A^{q}_2 
\bigg(\frac{l^2}{A_2}+\frac{m^2}{A_1}\bigg)^{-1}
\exp\bigg\{-\frac{A_1+A_2}{2\hb}\bigg\}\,dE,
$$
where $A_1=\frac{E}{2l}-m\alpha\zs{E}(x)$, 
$A_2=\frac{E}{2m}+l\alpha\zs{E}(x)$, $\alpha\zs{E}$ is taken from
{\rm(2.10)}, and $x=|z(a)|^2$. 
\end{lemma}

These are first steps to assign some geometry 
to the quantum algebra (2.4) and its irreducible
representations. 
The next step is to consider the multiplication operation
in this algebra. 

Note that linear operators in $\cP_{r,q,p}$ can be presented 
by their kernels. So, the algebra of operators is naturally
isomorphic to $\cS_{r,q,p}\od\cP_{r,q,p}\otimes\ol{\cP}_{r,q,p}$. 
The operator product is presented by the convolution of kernels
which is generated by pairing between $\ol{\cP}_{r,q,p}$ and
$\cP_{r,q,p}$ given by the scalar product (3.13).

The algebra $\cS_{r,q,p}$ consists of functions in $\ol{z}$,
$z$, they are polynomials on $\bR^2$. On this function space 
we have a noncommutative product (convolution),
but the unity element of this convolution is presented by the
function 
\begin{equation}
K=\sum^{r}_{k=0}\varphi^{(k)}\otimes \ol{\varphi^{(k)}},
\tag{4.3}
\end{equation}
where $\varphi^{(k)}$ is the orthonormal basis in $\cP_{r,q,p}$.
This function is called a reproducing kernel \cite{j38,j39}, 
it is independent of the choice of the basis $\{\varphi^{(k)}\}$.
From (3.15) we see the explicit formula for the reproducing 
kernel 
\begin{equation}
K=k(|z|^2),\qquad
k(x)\od\sum^{r}_{n=0}\hb^{(m-l)n}
\frac{q!(p+rm)!}{(q+nl)!(p+(r-n)m)!}x^n.
\tag{4.4}
\end{equation}

In order to give a Gelfand type spectral--geometric
interpretation of some algebra, we, first of all, have to ensure 
that the unity element of this algebra is presented by the unity
function.   
It is not so for the algebra $\cS_{r,q,p}$.
That is why we have to divide the ``kernel elements'' 
from $\cS_{r,q,p}$ by the reproducing kernel (4.4). 
The correct function algebra consists of ratios of the type 
\begin{equation}
f=\frac{\varphi\otimes\ol{\varphi'}}{K},
\tag{4.5}
\end{equation}
where $\varphi,\varphi'\in\cP_{r,q,p}$. The product of two
functions of this type generated by the convolution of kernels
is given by 
\begin{equation}
(f_1*f_2)(a)=\frac1{2\pi\hb}\int_{\text{phase space}} 
f^{\#}_1(a|b) f^{\#}_2(b|a) p_a(b)\,dm(b).
\tag{4.6}
\end{equation}

Here 
\begin{align}
dm(b)&\od L(b) K(b)\,db,
\tag{4.7}\\
p_a(b)&\od |K^{\#}(a|b)|^2 K(a)^{-1}K(b)^{-1},
\tag{4.8}
\end{align}
and the operation $f\to f^{\#}$ denotes the analytic continuation
holomorphic with respect to the ``right'' argument and
anti-holomorphic with respect to the ``left'' argument in the
notation $f^{\#}(\cdot|\cdot)$. 
The product (4.6) possesses the desirable property: 
$1*f=f*1=f$.

Let us look at formula (4.5). Since $f$ is going to be a
function on an invariant geometric space, $\varphi$ and
$\varphi'$ have to be sections of a Hermitian line bundle with
the curvature form 
\begin{equation}
\omega=i\hb\opa \pa \ln K\equiv igd\ol{z}\wedge dz.
\tag{4.9}
\end{equation}
Here $\pa$ denotes the differential by $z$.
Formula (4.9)  means that the measure $dm$ (4.7) is the
reproducing measure with respect to the K\"ahlerian form
$\omega$ in the sense \cite{j33}. 

Note that formula (4.9) defines both the {\it quantum form} 
$\omega$ and the {\it quantum metric} $g=g(|z|^2)$, 
$g(x)=\hb\frac{d}{dx}(x\frac{d}{dx}(\ln k(x)))$ via the
polynomial (4.4).

After the quantum form $\omega$ appears, the ``probability''
factor $p_a$ in the noncommutative product (4.6) can be written
as 
\begin{equation}
p_a(b)=\exp\bigg\{\frac i\hb\int_{\sum(a,b)}\omega\bigg\}.
\tag{4.10}
\end{equation}
Here $\sum(a,b)$ is a membrane in the complexified space
whose boundary consists of four paths connecting points 
$a \to b|a \to b \to a|b \to a$ along leaves of the complex
polarization and its conjugate \cite{j26,j40}.  

Note that the set of functions (4.10) makes up a resolution 
of unity:
\begin{equation}
\frac1{2\pi\hb}\int_{\text{phase space}}p_a\,dm(a)=1,
\tag{4.11}
\end{equation}
and each $p_a$ is the ``eigenfunction'' of the operators of left or
right multiplication:
\begin{equation}
f*p_a=f(\cdot|a)p_a,\qquad p_a*f=f(a|\cdot)p_a.
\tag{4.12}
\end{equation}
The details about such a way to establish a correspondence
between quantum algebras and K\"ahlerian geometry 
can be found in \cite{j33}. 

Let us discuss global aspects of this quantum geometry. 
The K\"ahlerian form $\omega$ (4.9) is actually well defined on
the compactified plane $\bR^2\cup\{\infty\}$ which includes the
infinity point $z=\infty$. 
To see this, we just can make the change of variables 
$z'=1/z$ and observe that $\omega$ is smooth near $z'=0$.

Thus the actual phase space is diffeomorphic to $\bS^2$ and we have 
\begin{equation}
\frac1{2\pi\hb}\int_{\bS^2}\omega=r,\qquad
\frac1{2\pi\hb}\int_{\bS^2}dm=r+1.
\tag{4.13}
\end{equation}

The first formula (4.13) follows from the fact that 
$K\sim\const\cdot|z|^{2r}$ as $z\to\infty$
(see in (4.4)).
It means that the cohomology class $\frac1{2\pi\hb}[\omega]$ is
integer, and this is the necessary condition for the Hermitian
bundle with the curvature $i\omega$ over $\bS^2$ to have global
sections \cite{j41}. 

The second formula (4.13) follows from the definition (4.3)
which implies 
$\frac1{2\pi\hb}\int dm=\sum^{r}_{k=0}\|\varphi^{(k)}\|^2$,
where the norm of each $\varphi^{(k)}$ is taken in the sense
(4.3) and is equal to $1$ by definition. 
The number $r+1$ in (4.13) is the dimension of the irreducible
representation of the resonance algebra.

We stress that the quantum K\"ahlerian form $\omega$, 
given by (4.4), (4.9), and the quantum measure $dm$, 
given by (4.7) and Lemma~4.1, are essentially different from
the classical form $\omega_0$ (2.9) and the classical Liouville
measure $dm_0=|\omega_0|$. 
The main difference is that $\omega$ is smooth 
and $dm$ is regular at poles while $\omega_0$ and $dm_0$ are
not. 
Some information regarding asymptotics of the quantum objects as
$\hb\to0$ and asymptotics near the poles 
is summarized in the following lemma. 

\begin{lemma}
{\rm(a)} 
In the classical limit $\hb\to0$, $E_{r,q,p}\to E>0$, 
out of neighborhoods of the poles $z=0$ and $z=\infty$ on the
sphere,  
the quantum geometrical objects are approximated by the
classical ones:
$$
\omega=\omega_0+O(\hb),\qquad  dm=dm_0(1+O(\hb)).
$$

{\rm(b)} 
The behavior of the quantum reproducing measure near the poles
is the following: 
\begin{equation}
\begin{aligned}
dm&\sim\const\cdot \frac{dx\wedge d\varphi}{x^{1-(q+1)/l}}
\qquad\text{as}\quad x\to0,\\
dm&\sim\const\cdot \frac{dx\wedge d\varphi}{x^{1+(p+1)/m}}
\qquad\text{as}\quad x\to\infty,
\end{aligned}
\tag{4.14}
\end{equation}
where $z=x^{1/2}\exp\{i\varphi\}$. 
Thus the reproducing measure has weak singularities at poles. 

{\rm(c)} Near the poles, the quantum K\"ahlerian form looks as 
\begin{align*}
\omega&\sim \hb^{m-l+1}\frac{(p+rm)!q!}{(p+rm-m)!(q+l)!}
id\ol{z}\wedge dz\qquad\text{as}\quad z\to0,
\\
\omega&\sim \hb^{l-m+1}\frac{p!(q+rl)!}{(p+m)!(q+rl-l)!}
\frac{id\ol{z}\wedge dz}{|z|^4}\qquad\text{as}\quad z\to\infty.
\end{align*}
Thus, near the poles, the asymptotics of $\omega$ as $\hb\to0$
is  
\begin{align}
\omega&\sim \const\hb^{1-l}id\ol{z}\wedge dz\qquad(z\sim 0),
\tag{4.15}\\
\omega&\sim \const\hb^{1-m}\frac{id\ol{z}\wedge dz}{|z|^4}
\qquad(z\sim\infty).
\nn
\end{align}
\end{lemma}

Comparing (4.14) with (2.11) we see that, near poles, 
$dm$ is not approximated by $dm_0$ as $\hb\to0$ if $q>0$ or
$p>0$. So, the usual deformation theory (starting with
classical data) 
{\it cannot be applied to compute the reproducing measure\/}
globally on the phase space.

Formulas (4.15) demonstrate that the quantum $\omega$
is not approximated by $\omega_0$ as $\hb\to0$
near the poles; the classical form $\omega_0$ must be
singular at $z=0$ if $l>1$ and be singular at $z=\infty$ if $m>1$. 
This statement is in  agreement with (2.11).

Note that the cohomology class of the classical symplectic form 
$\omega_0$ on the classical leaf with the quantized energy 
$E=E_{r,q,p}$ (3.12) is given by (2.7):
\begin{equation}
\frac1{2\pi\hb}\int_{\Omega_0}\omega_0
=r+\frac{q}{l}+\frac{p}{m}.
\tag{4.16}
\end{equation}
Here $r\sim \hb^{-1}$ is the main quantum number which controls
the dimension of the quantum Hilbert space $\cP_{r,q,p}$. 
The integers $q,p$ vary on the intervals (3.11), they control
the fine structure of the scalar product (4.2) in $\cP_{r,q,p}$.

In the case of ``ground states'', where $q=p=0$, 
the condition (4.16) becomes standard for the geometric
quantization. 
In the ``excited'' case where $q\geq1$ or $p\geq1$,
we observe something like an index contribution to the geometric
quantization picture appearing due to an additional holonomy
around the conical poles in $\Omega_0$. 
Because of these ``excitations'', the leaves $\Omega_0$ 
with quantized energies 
are distant from each other by $\frac{1}{l}$ or $\frac{1}{m}$  
fractions of the parameter~$\hb$.

To conclude this section, let us discuss what quantum leaves
of the algebra (2.4) are. To each element $\bF$ of the algebra
one can assign the corresponding operator $\bbf$ in the irreducible
representation. This operator acts in the Hilbert space $\cP_{r,q,p}$
of antiholomorphic sections over the phase space. 
Thus we can compose the function 
\begin{equation}
f\od\frac1K\bbf(K).
\tag{4.17}
\end{equation}
Here $K$ is the reproducing kernel (4.4) and
the operator $\bbf$ acts by $\ol{z}$. 
The function~$f$ (4.17) is called the {\it Wick symbol} of the
operator $\bbf$, for more details see in \cite{j3,j34,j42,j43}. 
The product of symbols in the sense of (4.6) corresponds to the
product of operators. Moreover, one can reconstruct the operator
by its symbol using the simple formula 
$$
\bbf=f(\overset{2}{\bz}{}^*,\overset{1}{\bz}),
$$
where $\bz^*$ is the operator of multiplication 
by $\ol{z}$ and~$\bz$ is the conjugate operator.

To generators of the algebra (2.4) we now can assign 
functions on the phase space:
\begin{equation}
a_j\od\frac1K\ba_j(K)\quad (j=1,2),
\qquad
a_\pm\od\frac1K\ba_\pm(K).
\tag{4.18}
\end{equation}
We can consider them as quantum analogs of the coordinate
functions $A_1$, $A_2$, $A_\pm=A_3\mp iA_4$ on classical
symplectic leaves of the Poisson algebra (2.5).

\begin{theorem}
{\rm(a)} The quantum coordinate functions obey the Casimir
identities 
\begin{gather*}
k a_1+m a_2=E_{r,q,p},
\\
a_+ * a_-=(a_1+\hb)*\dots*(a_1+m\hb)*a_2*(a_2-\hb)
*\dots*(a_2-l\hb+\hb).
\end{gather*}
Here $*$ is the quantum product {\rm(4.6)}.

{\rm(b)} In the classical limit $\hb\to0$ 
{\rm(}and $r\sim\hb^{-1}\to\infty${\rm)} the quantum coordinate
functions coincide with the classical coordinate functions
{\rm(2.12)} on the closure $\Omega$ {\rm(2.6)} of the symplectic
leaves~$\Omega_0$. 
\end{theorem}

Taking into account this theorem, we below identify the quantum
phase space $\bS^2$ with the closure $\Omega$ of the symplectic
leaf (2.6), where $E=E_{r,q,p}$.
We will call $\Omega$ endowed with this structure 
a {\it quantum leaf}.

Each element $\bF$ of the algebra (2.4) can be represented as a 
polynomial in generators:
\begin{equation}
\bF=F(\bA),\qquad 
\bA=(\overset{3}{\bA}_+,\overset{2}{\bA}_1,
\overset{2}{\bA}_2,\overset{1}{\bA}_-).
\tag{4.19}
\end{equation}
Here $F$ is a function on $\bR^4$.
The operation of multiplication of elements (4.19) determines 
a product operation $\odot$ in the algebra of polynomials
over~$\bR^4$:  
$$
F(\bA)G(\bA)=(F\odot G)(\bA)
$$
(see details in \cite{j32}).

Following \cite{j32,j33}, one can define the {\it quantum
restriction} of the function $F$ onto the leaves $\Omega$:
\begin{equation}
F\Big|_{\hat\Omega}\od \frac1K F(\ba)(K).
\tag{4.20}
\end{equation}

From \cite{j33} one known the following assertion.

\begin{theorem}
{\rm(a)} The quantum restriction {\rm(4.20)} 
$F\to F\Big|_{\hat\Omega}$ is a homomorphism of algebras{\rm:}
$$
(F\odot G)\Big|_{\hat\Omega}
=F\Big|_{\hat\Omega} * G\Big|_{\hat\Omega}.
$$

The equivalent formula for the  quantum restriction is 
$$
F\Big|_{\hat\Omega} =F(a*)1,
$$
where $a*$ are the operators of left multiplication by the
quantum coordinate functions $a=(a_+,a_1,a_2,a_-)$ {\rm(4.18)} 
in the algebra {\rm(4.6)}.

{\rm(b)} The asymptotics as $\hb\to0$ of the quantum restriction
can be derived from 
\begin{equation}
F\Big|_{\hat\Omega} =F\big(a-i\hb \ad_-(a)+O(\hb^2)\big)1
=F(a)+\hb e_1(F)+O(\hb^2).
\tag{4.21}
\end{equation}
Here $\ad_-(\cdot)$ denotes the anti-holomorphic part of 
the Hamiltonian field: $\ad_-(\cdot)=ig^{-1}\pa(\cdot)\opa$, 
where $g$ is the quantum metric {\rm(4.9)}.
The $\hb$-correction $e_1$ in {\rm(4.21)} is the second order
operator 
$e_1=\frac12\langle R\frac{\pa}{\pa a},\frac{\pa}{\pa a}\rangle$ 
determined by the symmetric tensor 
$R_{jl}=\Re(g^{-1}\pa a_j\opa a_l)$.
\end{theorem}

\section{Coherent states and gyron spectrum}

In the Hilbert space $\cP_{r,q,p}$ of anti-holomorphic sections
of the Hermitian line bundle with the curvature $i\omega$ over
the phase space $\Omega\approx\bS^2$ we have the irreducible
representation of the resonance algebra (2.4) by differential
operators 
\begin{align}
\ba_1&=\hb(rm+p)-\hb m\ol{z}\opa,\qquad 
\ba_2=\hb q+\hb l\ol{z}\opa,
\tag{5.1}\\
\ba_+&=\hb^m\prod^{m}_{j=1}(rm+p+j-m\ol{z}\opa)\cdot\ol{z},
\qquad
\ba_-=\frac{\hb^l}{\ol{z}}\prod^{l}_{s=1}(q-s+1+l\ol{z}\opa),
\nn
\end{align}
where $\opa=\pa/\pa\ol{z}$.

The unity section $1=\ol{z}^0$ is the vacuum vector for this
representation in the sense that it is the eigenvector 
of the operators $\ba_1$, $\ba_2$ and it is annulled by the
operator $\ba_-$. Now let us take the vacuum vector $\fP_0$
in the original Hilbert space $\cL=L^2(\bR^2)$ which correspond
to the representation (2.2):
$$
\bA_1\fP_0=\hb(rm+p)\cdot\fP_0,\qquad 
\bA_2\fP_0=\hb q\cdot\fP_0,\qquad
\bA_-\fP_0=0.
$$

\begin{definition}
{\rm
The {\it coherent states} of the algebra (2.4) is the
holomorphic family of vectors $\fP_z\in\cL$ defined by
$$
\fP_z=\sum^r_{n=0}\frac{q!}{(q+ln)!}
\bigg(\frac{z}{\hb^l}\bigg)^n\bA^n_+\fP_0,\qquad
z\in\bC.
$$}
\end{definition}

For each $a\in\Omega$ let us denote by $\boldsymbol\Pi_a$ the
projection onto the one-dimensional subspace in $\cL$
generated by $\fP_{z(a)}$. 
We call $\Pi_a$ a {\it coherent projection}.

Regarding these definitions, may be, it is useful to note the
following:
if one takes the Hilbert space $\cP_{r,q,p}$ instead of $\cL$
and the vacuum~$1$ instead of $\fP_0$, 
then instead of coherent states $\fP_z$ and
the coherent projection $\boldsymbol\Pi_a$ 
one would see the reproducing kernel $K^{\#}(\cdot|z)$ 
and the probability function~$p_a$.

In the following theorem we collect the basic properties of the
coherent states $\fP_z$. In the general context of quantization
theory, see more details in \cite{j32,j33,j34}.

\begin{theorem}
{\rm(a)} The scalar product of two coherent states coincides with
the reproducing kernel {\rm(4.4):}
$$
\|\fP_{z(a)}\|^2=K(a),\qquad a\in\Omega.
$$

{\rm(b)} One has the resolution of unity by coherent
projections{\rm:} 
$$
\frac1{2\pi\hb}\int_{\Omega}\boldsymbol\Pi_a\,dm(a)=\bI_{r,q,p}.
$$
Here $\bI_{r,q,p}$ is the projection in $\cL$ onto the Hilbert
subspace $\cL_{r,q,p}$ spanned by all vectors
$\bA^n_+\fP_0$, $n=0,\dots,r$. 

{\rm(c)} The whole Hilbert space $\cL$ is the direct sum 
of the irreducible subspaces:
$$
\cL=\bigoplus_{\substack{r\geq0\\ 0\leq q\leq l-1\\ 0\leq p\leq m-1}}
\cL_{r,q,p}.
$$

{\rm(d)} The coherent transform
$\cL_{r,q,p}\overset{\nu}{\to}\cP_{r,q,p}$ defined by  
\begin{equation}
\nu(\psi)(\ol{z})=(\psi,\fP_z),
\tag{5.2}
\end{equation}
has the inverse
\begin{equation}
\nu^{-1}(\varphi)=\frac1{2\pi\hb}\int_\Omega
\frac{\fP\otimes \varphi}{K}\,dm.
\tag{5.3}
\end{equation}
The mappings {\rm(5.2), (5.3)} intertwine the representations 
{\rm(2.2)} and {\rm(5.1)} of the algebra {\rm(2.4)}.

{\rm(e)} Let $\bF$ be an element of the algebra {\rm(2.4)} 
realized in the Hilbert space $\cL$ via the generators {\rm(2.2)}
as in {\rm(4.19)}, and let $\bbf=\nu\circ \bF\circ \nu^{-1}$ 
be the coherent transformation of $\bF$ realized in the Hilbert
space $\cP_{r,q,p}$. 
Then the Wick symbol $f$ {\rm(4.17)} coincides with the Wick
symbol of $\bF$ given by 
$$
f(a)=\tr(\bF\boldsymbol\Pi_a),\qquad a\in\Omega.
$$
The operators $\bF,\bbf$ are reconstructed via their symbols using 
the formulas 
\begin{equation}
\bF=F(\bA)=f(\overset{2}{\bz}{}^*,\overset{1}{\bz}),
\qquad 
\bbf=F(\ba)=f(\overset{2}{\ol{z}},\overset{1}{\ol{z}}{}^*),
\tag{5.4}
\end{equation}
where $\bz$ is the operator of complex structure {\rm(4.1)},
$\ba$ are the operators of irreducible representation {\rm(5.1)}.
The Wick symbol of the coherent projection $\boldsymbol\Pi_a$ is
the probability function $p_a$ {\rm(4.10)}.
\end{theorem}

Now following \cite{j36},\cite{j44}--\cite{j48} 
we explain how to reduce the coherent transform to closed curves
(Lagrangian submanifolds) in the phase space.

Let $\Lambda\subset\Omega$ be a smooth closed curve,
which obeys the quantization condition 
\begin{equation}
\frac1{2\pi\hb}\int_{\Sigma}\Big(\omega-\frac{\hb}{2}\rho\Big)
-\frac{1}{2}\in\bZ,
\tag{5.5}
\end{equation}
where $\omega=igd\ol{z}\wedge dz$ is the quantum K\"ahlerian
form (4.9), $\rho=i\opa\pa \ln g$ is the quantum Ricci form,
and $\Sigma$ is a membrane in $\Omega$ with the boundary 
$\pa\Sigma=\Lambda$. 

We choose certain parameterization of the curve
expressed via the complex coordinate on the leaf as follows:
$$
\Lambda=\{z=z(t)\mid 0\leq t\leq T\},
$$
and define the following basis of smooth functions on the curve:
\begin{equation}
\phi^{(j)}(t)=\sqrt{\ol{\dot z(t)}}
\exp\bigg\{-\frac{i}{\hb}\int^{t}_{0}
\Big(\ol{\theta}-\frac{\hb}{2}\ol{\varkappa}\Big)\bigg\} 
\varphi^{(j)}(\ol{z}(t)),\qquad j=0,\dots,r.
\tag{5.6}
\end{equation}
Here $\theta\od i\hb\pa\ln K$ and $\varkappa=i\pa\ln g$ are
primitives of the quantum K\"ahlerian form $\omega=d\theta$
and the quantum Ricci form $\rho=d\varkappa$, 
the integral in (5.6) is taken over a segment 
of the curve $\Lambda$, and the monomials
$\varphi^{(j)}$ are defined in (3.15). 

Let us denote by $\cL_\Lambda$ the vector subspace in
$C^\infty(\Lambda)$ spanned by $\phi^{(j)}$ ($j=0,\dots,r$)  
and introduce the Hilbert structure in $\cL_\Lambda$ by means of
the following norm:
\begin{equation}
\|\phi\|\zs\Lambda\od\frac1{\root{4}\of{2\pi\hb}}
\bigg(\sum^{r}_{j=0}
\big|(\phi,\phi^{(j)})_{L^2}\big|^2\bigg)^{1/2},
\tag{5.7}
\end{equation}
where the scalar product $(\cdot,\cdot)_{L^2}$ is taken
in the $L^2$-space over $\Lambda$.

For any smooth function $\phi\in C^\infty(\Lambda)$ 
we define  
\begin{equation}
\mu\zs\Lambda(\phi)=\frac1{\root{4}\of{2\pi\hb}}
\int\zs\Lambda\phi(t)
\sqrt{\dot z(t)}
\exp\bigg\{\frac{i}{\hb}\int^{t}_{0}
\Big(\theta-\frac{\hb}{2}\varkappa\Big)\bigg\}
\fP_{z(t)}\,dt,
\tag{5.8}
\end{equation}
where $\fP\in\cL$ are coherent states of algebra (2.4)
corresponding to its $(r,q,p)$-irreducible representation.

\begin{theorem}
{\rm(a)}
The mapping $\mu\zs\Lambda$ defined by {\rm(5.8)} 
is an isomorphism of Hilbert spaces
$$
\mu\zs\Lambda:\,\cL_\Lambda\to \cL_{r,q,p}\subset \cL.
$$ 

{\rm(b)}
Under the isomorphism {\rm(5.8)} the representation of the
algebra {\rm(2.4)} in the Hilbert space $\cL$ is transformed to
the irreducible representation in the Hilbert space
$\cL_\Lambda${\rm:} 
\begin{equation}
\bF\to\bF_\Lambda\od 
\mu^{-1}\zs\Lambda\circ\bF\circ\mu\zs\Lambda.
\tag{5.9}
\end{equation}

{\rm(c)}
In the classical limit as $\hb\to0$ the Hilbert structure
{\rm(5.7)} coincides with the $L^2$-structure{\rm:}
\begin{equation}
\|\phi\|\zs\Lambda=\bigg(\int_\Lambda |\phi(t)|^2\,dt\bigg)^{1/2}+O(\hb).
\tag{5.10}
\end{equation}

{\rm(d)}
Let $f$ be the Wick symbol {\rm(5.4)} of the operator $\bF$,
then the asymptotics of the operator {\rm(5.9)} as $\hb\to0$ is
given by 
\begin{equation}
\bF_\Lambda=\cF\Big|_{\Lambda}
-i\hb\Big(v+\frac{1}{2}\div v\Big)+O(\hb^2).
\tag{5.11}
\end{equation}
Here $\cF=f-\frac{\hb}{4}\Delta f$, 
by $\Delta$ we denote the Laplace operator 
with respect to the quantum K\"ahlerian metric $g$,
and 
$v=\ad_+(\cF)\big|\zs\Lambda$ is the restriction to
$\Lambda$ of the holomorphic part of the  Hamiltonian field  
$\ad_+(\cF)=-ig^{-1}\opa\cF\cdot\pa$.
\end{theorem}

The next terms of the asymptotic expansion (5.11) are also known
(see in \cite{j36}).

In Theorem~5.2, the curve $\Lambda$ is arbitrary except it has
to obey the quantization condition (5.5).

Let us now choose $\Lambda$ specifically to be a closed curve 
on the energy level
\begin{equation}
\Lambda\subset\{\cF=\lambda\},
\tag{5.12}
\end{equation}
and choose the coordinate $t$ to be time 
on the trajectory $\Lambda$ of the Hamiltonian field $\ad(\cF)$.
Then $v=\ad(\cF)\big|\zs\Lambda=\frac{d}{dt}$, 
$\div v=0$, and we have 
\begin{equation}
\bF_\Lambda=\lambda-i\hb\frac{d}{dt}+O(\hb^2).
\tag{5.13}
\end{equation}
This formula implies the asymptotics of eigenvalues 
of the operator $\bF_\Lambda$:
\begin{equation}
\lambda+\hb\frac{2\pi k}{T}+O(\hb^2),
\tag{5.14}
\end{equation}
where $T=T(\lambda)$ is the period of the trajectory 
$\Lambda=\Lambda(\lambda)$ (5.12) and $\lambda$ is determined by
the quantization condition (5.5). 

Note that the contribution $\frac{2\pi k}{T}$ added to $\lambda$
in (5.14) can be transformed to adding the number $k$
to the integer number on the right-hand side of condition (5.5).
Thus one can omit the summand $\hb\frac{2\pi k}{T}$ in (5.16) 
without loss of generality.

\begin{corollary}
Let $\bF$ be an operator commuting with the oscillator $\bE$ 
{\rm(2.1)}.
Up to $O(\hb^2)$, the asymptotics of its eigenvalues $\lambda$ 
is determined by the quantization condition{\rm:}
\begin{equation}
\frac1{2\pi\hb}\int_{\Sigma}\Big(\omega-\frac{\hb}{2}\rho\Big)
-\frac12\in\bZ.
\tag{5.15}
\end{equation}
Here $\Sigma$ is a membrane in $\Omega$ with the boundary 
$\Lambda=\pa\Sigma$ {\rm(5.12);} 
the curve $\Lambda$ is the energy level of the function 
$\cF=f-\frac{\hb}{4}\Delta f$, where $f$ is the Wick
symbol of $\bF$ and $\Delta$ is the Laplace operator. 
The operator  $\Delta$ and the forms $\omega$, $\rho$
are generated by the quantum K\"ahlerian metric $g$ {\rm(4.9)}. 
\end{corollary}

Now we can apply the obtained results in studying quantum
gyrons. Let one has the Hamiltonian of the type 
\begin{equation}
\bE+\ve\bB,
\tag{5.16}
\end{equation}
where $\bE$ is the oscillator (2.1) and $\bB$ is a perturbation 
presented as a function in operators $\bb, \bb^*$,
\begin{equation}
\bB=\sum\beta_{\mu,\nu}{\bb^*}^\nu \bb^\mu.
\tag{5.17}
\end{equation}
There is an {\it operator averaging}  procedure 
\cite{j47,j48}, which is a unitary
transformation reducing (5.16) (up to $O(\ve^N)$) to the
Hamiltonian  
\begin{equation}
\bE+\ve\bB\sim\bE+\ve\bF_N+O(\ve^N),\qquad
[\bF_N,\bE]=0.
\tag{5.18}
\end{equation}
For instance, if $N=1$, then 
\begin{equation}
\bF_1=\sum_{l\nu_1+m\nu_2=l\mu_1+m\mu_2}
\beta_{\mu,\nu}{\bb^*}^{\nu}\bb^\mu
\tag{5.19}
\end{equation}
(see also the Appendix in \cite{j25}).
For any $N\ge1$ in (5.18), the operator $\bF_N$, 
commuting with $\bE$, 
is uniquely determined 
and can be presented in the form (4.19):
$$
\bF_N=F_N(\bA),
$$
and after this in the form (5.4):
\begin{equation}
\nu\circ \bF_N\circ \nu^{-1}= F_N(\ba)
=f\zs N(\overset{2}{\ol{z}},\overset{1}{\ol{z}}{}^*).
\tag{5.20}
\end{equation}

Thus the study of the operator (5.16) up to $O(\ve^N)$ is
reduced to the study of the properties of the operator (5.20) in
each irreducible representation of the algebra (2.4).

The symbols $F_N$ or $f\zs N$ are {\it gyron Hamiltonians}. 
In the $(r,q,p)$-irreducible representation, the gyron is
described by the operator 
$F_N(\ba)=F_N(\overset{3}{\ba}_+,\overset{2}{\ba}_1,
\overset{2}{\ba}_2,\overset{1}{\ba}_-)$ acting in $\cP_{r,q,p}$, 
where the generators $\ba$ are given by (5.1).

In the semiclassical approximation $\hb\to0$ 
the gyron system can be reduced to (5.11) and even to (5.13)
over the 
trajectory $\Lambda$ of the effective Hamiltonian 
$\cF_N=f\zs N-\frac{\hb}{4}\Delta f\zs N+O(\hb^2)$ 
on the leaf $\Omega\approx\bS^2$. The asymptotics of the
gyron spectrum was described in Corollary~5.3
by means of the membrane versions (5.15) 
of the Bohr--Sommerfeld quantization condition.

The quantum K\"ahlerian geometry (via the measure $dm$ and 
the forms $\omega,\rho$) 
is essentially presented in all these results regarding 
the gyron spectrum. 

The gyron is a model. It is very simple, since it arises from 
the ``textbook'' oscillator Hamiltonian. At the same time, it
already contains many nontrivial aspects of the quantization
theory and, of course, it has a variety of important physical
applications. About more complicated models of this type and
about further ideas on the quantum geometry we refer to 
\cite{j25,j32,j33}, \cite{j51}--\cite{j63}.

{\bf Acknowledgements.}
The author is grateful to V.~P.~Maslov and E.~M.~No\-vi\-ko\-va for
very useful discussions and help.

\end{document}